\documentclass[a4paper,10pt,reqno]{amsart}        
\usepackage{verbatim}       
\usepackage{amssymb}       
\usepackage{mathrsfs}       
       
\usepackage{enumerate}       
\usepackage[active]{srcltx}       
\numberwithin{equation}{section}

\usepackage{tikz}
\usetikzlibrary{positioning}
\usepackage{marvosym}
\usetikzlibrary{arrows.meta,calc}

 \makeatletter      
  \@namedef{subjclassname@2020}{%
  \textup{2020} Mathematics Subject Classification}    
   \makeatother    
       
\usepackage{t1enc}       
\usepackage[utf8x]{inputenc}

\newtheorem{theorem}{Theorem}[section]        
\newtheorem{lemma}[theorem]{Lemma}       
\newtheorem{proposition}[theorem]{Proposition}

\theoremstyle{definition}       
\newtheorem{definition}[theorem]{Definition}       
\theoremstyle{remark}       
         
\newcommand{\mc}[1]{\mathcal{#1}}       
       
\newcommand{\mbb}[1]{\mathbb{#1}}

\newcommand{\setm}{\setminus}       
\newcommand{\empt}{\emptyset}       
\newcommand{\subs}{\subset}

\def\<{\left\langle}       
\def\>{\right\rangle}       
\author[I. Juh\'asz]{Istv\'an Juh\'asz}       
\address       
      { Alfr{\'e}d R{\'e}nyi Institute of Mathematics, Hun-Ren}
       
\email{juhasz@renyi.hu}

\author[L. Soukup]{Lajos Soukup}       
\address       
      { Alfr{\'e}d R{\'e}nyi Institute of Mathematics,  Hun-Ren
}       
\email{soukup@renyi.hu}       
       
\author[Z. Szentmikl\'ossy]{Zolt\'an Szentmikl\'ossy}       
\address{Eötvös University of Budapest}       
\email{szentmiklossyz@gmail.com}       
       

\title[On the companion of \DRS\ spaces]{On the companion of spaces having dense, relatively countable compact subspaces }       
\thanks{The preparation of this paper was
supported by  OTKA grant   K129211}       
\date{\today}

\newcommand{\DRC}{{DRC}}
\newcommand{\DRCo}{{DRC${}_{\omega}$}}
\newcommand{\DRS}{{DRS}}
\newcommand{\DRSo}{DRS${}_{\omega}$}

\newcommand{\descri}[1]{\<X^{#1},\{B^{#1}_i:i\in{\nu}^{#1}\}\>}
\newcommand{\descris}[1]{\<X^{#1},\mc B^{#1}\>}

\newcommand{\appros}[1]{\<X^{#1},\mc B^{#1},  \mc F^{#1}\>}

\newcommand{\dset}{\operatorname{\mc D}}

\begin{document}     
\begin{abstract}
A topological space is said to be \DRC\ (\DRS)\ iff it possesses a dense, 
relatively countably compact (or relatively sequentially compact, respectively) subspace. 

   The concept of 
 selectively pseudocompact game Sp(X) and the 
 selectively sequentially pseudocompact game Ssp(X) were
 introduced by   Dorantes-Aldama and Shakhmatov. 
 They explored the relationship between the existence of a winning strategy
 and a stationary winning strategy for player P in these games. 
In particular,  they observed that there exists  a stationary winning strategy
 in the game Sp(X) (Ssp(X)) for Player P iff $X$ is \DRC\ (or \DRS, respectively).

In this paper we introduce natural weakening of the properties 
\DRC\ and \DRS: 
a space $X$ is {\em \DRCo\ (\DRSo)  } iff there is a 
 sequence $\<D_n:n\in {\omega}\>$ of dense subsets of $X$ such that every sequence 
$\<d_n:n\in {\omega}\>$ with 
 $d_n\in D_n$ has an accumulation point (or contains a convergent subsequence, respectively).

These properties  are also equivalent 
to the existence of  some limited knowledge winning strategy on the 
corresponding games $Sp(X)$ and $Ssp(X)$.

Clearly, \DRS\ implies \DRC\ and \DRSo, \DRC\ or \DRSo\ imply \DRCo.
The main part of this paper is devoted to prove  that apart from these trivial implications, consistently there are no 
other  implications  between these properties. 
\end{abstract}

\maketitle

\section{Introduction}

The notion of {\em pseudocompactness} was introduced by Hewitt in \cite{He48}.
The concept of relatively countably compact subspaces were explored by 
Marjanovic  in \cite{Ma71} to show that a $\Psi$-space is pseudocompact. 
Berner \cite{Be81} constructed pseudocompact spaces with and without having dense,   
relatively countably compact subspaces. 

To simplify the formulation of our results, a topological space is said to be \DRC\ (\DRS)\ iff it possesses a dense, 
relatively countably compact (or relatively sequentially compact, respectively) subspace.

The concept of selective pseudocompact and selectively sequentially pseudocompact  
spaces were introduced in \cite{DoSh17}. In a subsequent work, \cite{DoSh20},
the same authors introduced 
two related topological games and explored  
the connection between these classes of spaces and the existence of certain type of winning strategies in the defined
games.    

To start with, we recall some definitions from \cite{DoSh17} and \cite{DoSh20}.
Given a topological space $X$ and a sequence $\vec a\in {}^{\omega}X$ write 
\begin{displaymath}
acc(\vec a)=\{x\in X:\{n\in {\omega}:\vec a(n)\in U\}\text{ is infinite 
for each open }U\ni x\}.
\end{displaymath}

A subspace $Y\subs X$   is {\em relatively countably compact} 
iff  $acc(\vec y)\ne \empt $ for   every $\vec y\in {}^{\omega}Y$.  
 A subspace $Z\subs X$ is {\em relatively sequentially  compact} 
 iff every sequence  $\vec z\in {}^{\omega}Z$ 
 contains   a subsequence converging to some point in $X$. 

 A space $X$ is {\em selectively pseudocompact} 
 (\cite[Def 2.2.]{DoSh17})    if 
 for every $\vec U\in {}^{\omega}({\tau}^+_X)$ there is 
 $\vec x\in {}^{\omega}X$ with $\vec x(n)\in \vec U(n)$
 such that $acc(\vec x)\ne \empt$.

 A space $X$ is {\em selectively sequentially pseudocompact}
 (\cite[Def 2.3]{DoSh17}) iff
 for every $\vec U\in {}^{\omega}({\tau}^+_X)$ there is 
 $\vec x\in {}^{\omega}X$ with $\vec x(n)\in \vec U(n)$
 such that $\vec x$
 contains a converging subsequence.

The following games were introduced in \cite[Definition 5.1]{DoSh20}.
   Given a  space $X$ define the games $Sp(X)$ and $Ssp(X)$ between 
   players O and P as follows. The games are played in ${\omega}$ turns.
   In both game in the $n^{th}$ turn O picks a non-empty open set $U_n$, and 
   then P picks a point $x_n\in U_n$.
         
   In the game $Sp(X)$ Player P wins   iff  
 $acc(\<x_n:n\in {\omega}\>)\ne \empt$.

In the game $Ssp(X)$ player P wins iff the sequence 
$\<x_n:n\in {\omega}\>$  contains  a convergent subsequence.

In \cite[Theorem 4.5]{DoSh20} the authors observed  
that the statements  (a) and (b) below  are equivalent:
(a) {\em the  space $X$ is contains a dense, 
relatively countably compact subset,} (b) {\em P has a stationary 
winning strategy in the game Sp(X)}. Let us recall that strategy of P
is {\em stationary} iff the moves of P depends on only from the last move
of the opponent.  

They also showed  that  (c) and (d) below  are similarly  equivalent:
(c) {\em the  space $X$  contains a dense, 
relatively sequentially  compact subset,} (b) {\em P has a stationary 
winning strategy in the game Spp(X)}.

 In \cite{DoSh20} the authors refrain from introducing specific names 
for spaces having property (a) or property (c). 
To simplify the formulation of our results we introduce the following names 
and notations. We say that 
a space $X$ is {\DRC} (``{\em densely relatively countably compact}'')  
iff it contains a dense, relatively countably compact subset,
and  $Y$ is {\DRS} (``{\em densely relatively sequentially compact}'') 
iff it contains a dense, relatively sequentially compact subset. 

We write $P\uparrow_{1}Sp(X)$  iff P has a stationary winning strategy in the
game $Sp(X)$. We define  $P\uparrow_{1}Spp(X)$ similarly.

Next, we introduce 
 natural weakenings of the \DRC\ and \DRS\ properties  which are still stronger than  the 
selectively pseudocompactness and selectively sequentially pseudocompactness, respectively.  

\begin{definition}\label{df:drco}
   Let $X$ be a topological space.
 \begin{enumerate}[(i)]
 \item $X$ is {\em \DRCo} iff there is a 
 sequence $\<D_n:n\in {\omega}\>$ of dense subsets of $X$ such that every sequence 
$\<d_n:n\in {\omega}\>$ with 
 $d_n\in D_n$ has an accumulation point.
 \item $X$ is {\em \DRSo} iff there is a 
 sequence $\<d_n:n\in {\omega}\>$ of dense subsets of $X$ such that 
 every sequence 
$\<d_n:n\in {\omega}\>$ with 
 $d_n\in D_n$ contains a convergent subsequence.
\end{enumerate}  
\end{definition}

As it turns out, these properties have characterizations using the 
games Sp and Spp. To formulate our observation
we need to introduce the following  types of 
strategies of games which use only restricted information.

Write $P\uparrow_{1,n}Sp(X)$   iff P has a  winning strategy in the  game such
$Sp(X)$  that  the $n^{th}$ move of $P$
depends on only  $n$ and   $U_n$. 
We define  $P\uparrow_{1,n}Ssp(X)$ similarly.

\begin{proposition}
A topological space $X$ is \DRCo\ iff $P\uparrow_{1,n}Sp(X)$.
A topological space $Y$ is \DRCo\ iff $P\uparrow_{1,n}Sp(Y)$.
\end{proposition}

\begin{proof}
Assume that $X$ is \DRCo\ witnessed by the sequence $\<D_n:n\in {\omega}\>$ of dense subsets of 
$X$. Let $P$ play the following strategy: in the $n^{th}$-turn if O pick $U_n$, then let P choose 
an arbitrary $d_n\in D_n\cap U_n$. Then the sequence $\<d_n:n<{\omega}\>$ has an accumulation point so 
$P$ wins. Moreover, the move of P depends on only $n$ and $U_n$.

Assume know that ${\sigma}: {\omega}\times {\tau}^+_X\to X$ is a winning strategy of P.
For $n\in {\omega}$, let
\begin{displaymath}
D_n=\{{\sigma}(n,U):U\in {\tau}^+(X)\}.
\end{displaymath} 
Then $D_n$ is dense because ${\sigma}(n,U)\in U$.

Moreover, if $\<d_n:n\in {\omega}\>$ is a sequence with $d_n\in D_n$, then we can pick 
$U_n\in {\tau}^+_X$ with $d_n={\sigma}(n,U_n)$.
If O plays $U_n$ in the $n^{th}$ turn, then P produces the sequence $\<d_n:n<{\omega}\>$.
Since P wins, $\<d_n:n<{\omega}\>$ has accumulation point. So the sequence 
$\<D_n:n<{\omega}\>$ witnesses that $X$ is \DRCo.

The second equivalence can be proved similarly. 
\end{proof}

Figure  \ref{fig1} summarizes the equivalences  and the straightforward   implications  between these properties. 

\tikzset{tul/.style={draw,minimum width=2cm,font=\small},el/.style={->,-Stealth},
ella/.style={<->,Stealth-Stealth},
uj/.style={fill=gray!50}}

\begin{figure}[h]\label{fig1}

\begin{tikzpicture}[every node/.style={scale=0.8}]

\node[tul]  (drc)  {\DRC};
\node[tul, left =2cm of   drc] (drs)  {\DRS};

\node[tul, right =1cm of   drc] (stay)  {$P\uparrow_1Sp(X)$};
\node[tul, left =1cm of   drs] (sstay)  {$P\uparrow_1Ssp(X)$};

\node[tul,uj,,below = 2cm of  drc] (drco)  {\DRCo};
\node[tul,uj,below = 2cm of drs] (drso)  {\DRSo};

\node[tul,uj, below =2cm of   stay] (stayo)  {$P\uparrow_{1,n}Sp(X)$};
\node[tul,uj, below =2cm of   sstay] (sstayo)  {$P\uparrow_{1,n}Ssp(X)$};

\node[below =1cm of   sstay]  (egye) {};
\node[below =1cm of   stay]  (egyw) {};
\draw[dashed] (egye) -- (egyw);
\node at (egye) {1};

\node[above right =0.9cm of   drs]  (nulla) {};
\node[below right =0.9cm of   drso]  (nullaw) {};
\draw[dashed] (nulla) -- (nullaw);
\node[left= 0mm of nullaw]  {0};

\node[below right =0.7cm of   sstay]  (ketto) {};
\node[above right =0.7cm of   drs]  (kettow) {};
\draw[dashed] (ketto) -| (kettow);
\node[above= 0mm of ketto]  {2};

\node[above left =0.7cm of   stayo]  (harom) {};
\node[below left =0.7cm of   drco]  (haromw) {};
\draw[dashed] (harom) -| (haromw);
\node[below= 0mm of harom]  {3};

\draw (sstay) edge[ella] (drs) (drc) edge[ella] (stay);

\draw (sstayo) edge[ella] (drso) (drco) edge[ella] (stayo);

\draw  (drs) edge[el] (drso) ;

\draw    (drc) edge[el] (drco)  ;

\draw  (drs) edge[el] (drc) (drso) edge[el] (drco) ;

\end{tikzpicture}

\caption{\DRC\ and its companion}
\end{figure}
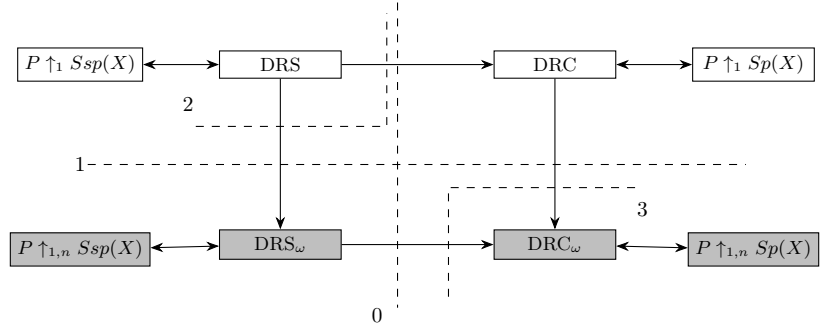

The aim of this paper is to show that these implications  
are not  reversible. 
Since there are countably compact spaces without having any convergent sequences 
(for example, ${\omega}^*$) , 
\DRC\ does not imply \DRSo.

The main result of this paper is Theorem  \ref{tm:mainmain}
below, which implies 
that the other implications are also not reversible.  

Given a topological property $Q$, we say that a space $X$ is {\em hereditary Q} iff every 
non-empty regular closed subset of $X$ has property $Q$.
For example, a pseudocompact space is ``hereditary pseudocompact'', although a closed subspace of a  
pseudocompact space is not necessarily pseudocompact.

\begin{theorem}[CH]\label{tm:mainmain}
   There is a crowded, 0-dimensional Hausdorff space $X$ of cardinality ${\omega}_1$, and 
   $X$ has three dense subspaces $Z_1$,
   $Z_2$ and $Z_3$ such that
   \begin{enumerate}[(1)]
   \item $Z_1$ is hereditary (\DRSo, but not \DRC),
   \item $Z_2$ is hereditary (\DRSo\ and \DRC, but not \DRS),
   \item $Z_3$ is hereditary (\DRCo, but neither \DRC, nor \DRSo).
   \end{enumerate} 
   \end{theorem}

   \begin{proof}[Proof of Theorem \ref{tm:mainmain}]
      Let $K,L,M$ be disjoint countable sets, let $\{K_n\}_{n\in {\omega}}$
      and $\{L_n\}_{n\in {\omega}}$ be partitions of $K$ and $L$ , 
      respectively, into infinite pieces.
   Write
   \begin{align*}
   \mc K&=\{K'\subs K:\forall n\  |K'\cap K_n|<{\omega}\},\\
   \mc L&=\{L'\subs L:\forall n\  |L'\cap L_n|<{\omega}\},\\
   \mc M&= {[M]}^{{\omega}}.
   \end{align*}

   The underlying set of our topological space  will be
   \begin{displaymath}
   K\cup L\cup M\cup {\omega}_1.
   \end{displaymath} 
   
   \begin{proposition}[CH]\label{lm:main}
      There exists  a 0-dimensional $T_2$ topological space $X=\<K\cup L\cup M\cup {\omega}_1,{\tau}\>$
      such that 
      \begin{enumerate}[(a)]
      \item \label{dense-level} $\forall U\in {\tau}^+$ 
      \begin{displaymath}
      |U\cap {\omega}_1|={\omega}_1\ \land |U\cap M|={\omega}\land   
      \forall^\infty n\ |U\cap K_n|=|U\cap L_n|={\omega},
      \end{displaymath}
      \item  \label{kconv} every $A\in \mc K$ contains a convergent subsequence. 
      \item \label{onlyconv} If $A\in {[X]}^{{\omega}}$ is convergent, then   
      $A\setm K$ is finite and $A\cap K\in \mc K$,
      \item \label{klacc} if $B\in \mc K\cup \mc L$, then $acc(B)=acc(B)\cap {\omega}_1\ne \empt$,
      \item \label{macc} if $B\in \mc M$, then $acc(B)\cap {\omega}_1\ne \empt$,
      \item \label{noacc}
$acc(K_n)=acc(L_n)=acc(\Gamma)=\empt$ for each $n<{\omega}$ and 
      $\Gamma\in {[{\omega}_1]}^{{\omega}}$  

   \end{enumerate}
      \end{proposition}
We show that  the  space $X$  from Proposition \ref{lm:main} by taking  
\begin{align*}
   Z_1=&K\cup{\omega}_1,\\
   Z_2=&(K\cup M)\cup {\omega}_1,\\
   Z_3=&L\cup {\omega}_1
   \end{align*}
 satisfies the requirements of Theorem \ref{tm:mainmain}.
 First, $X$ is 0-dimensional, $T_2$ and $|X|={\omega}_1$.

\smallskip 

\noindent 
 (1) Assume that $Y\subs Z_1$ is relatively countably compact. Then, by (\ref{noacc}), $Y\cap {\omega}_1$
is finite, and  $Y\cap K_n$ is also finite for each $n$. Thus, $Y\cap K\in \mc K$.
Hence, $acc(Y\cap K)\subs {\omega}_1$, and so  $Y\cap K$ is nowhere dense by (\ref{dense-level}).
Since $Y\cap {\omega}_1$ is finite, it follows that $Y$ is nowhere dense, and so  
no regular closed subset of $Z_1$ is \DRC.  

Next we show that $Z_1$ is \DRSo.
Fix first  a partition  $\{I_k:k\in {\omega}\}$ of ${\omega}$ into infinite pieces, and  
write $D_k=\bigcup\{K_n:n\in I_k\}$. Then, by (\ref{dense-level}),  every $D_k$ is dense.

If $d_k\in D_k$, then $\{d_k:k\in {\omega}\}\in \mc K$, so it contains a convergent subsequence by 
(\ref{kconv}). Thus,
the sequence $\{D_k:k\in {\omega}\}$ witnesses that $Z_1$ is \DRSo.

\smallskip

\noindent 
(2) Assume that $Y\subs Z_2$ is relatively sequentially compact.
Then,  $Y\setm K$ is finite, and $Y\cap K_n$ is finite for each $n$ by (\ref{onlyconv}).
Thus, $Y\cap K\in \mc K$. 
Hence $acc(Y\cap K)\subs {\omega}_1$ by (\ref{klacc}). Thus $Y\cap K$ is nowhere dense by 
(\ref{dense-level}).
Since $Y\setm K$ is finite, it follows that $Y$ is nowhere dense, and so 
no regular closed subset of $Z_2$ is \DRS.

The subspace $Z_2$ is \DRC, because $M$ is relatively countably compact by 
(\ref{klacc}) and $M$ is dense by (\ref{dense-level}).

Since $Z_1$ is a dense subspace of $Z_2$ because $K$ is dense in $X$ by (\ref{dense-level}),
and $Z_1$ is \DRSo
, it follows that $Z_2$ is also \DRSo.

\smallskip

\noindent 
(3) By (\ref{onlyconv}), $Z_3$ does not contain convergent sequences, so 
no subspace of  $Z_2$ is \DRSo.

Assume that $Y\subs Z_3$ is relatively countably compact. 
Then, by (\ref{noacc}), $Y\cap {\omega}_1$
is finite, and  $Y\cap L_n$ is also finite for each $n$. Thus, $Y\cap L\in \mc L$.
Hence $acc(Y\cap L)\subs {\omega}_1$, and so  $Y\cap L$ is nowhere dense by (\ref{klacc}).
Since $Y\cap {\omega}_1$ is finite, it follows that $Y$ is nowhere dense, and so  
no regular closed subset of $Z_3$ is \DRC. 

Next we show that $Z_3$ is \DRCo.
Fix first  a partition  $\{I_k:k\in {\omega}\}$ of ${\omega}$ into infinite pieces, and  
write $D_k=\bigcup\{K_n:n\in I_k\}$. Then, by (\ref{dense-level}),  every $D_k$ is dense.

If $d_k\in D_k$, then $\{d_k:k\in {\omega}\}\in \mc L$, so it 
has an accumulation point (\ref{klacc}). Thus,
the sequence $\{D_k:k\in {\omega}\}$ witnesses that $Z_3$ is \DRCo.

\bigskip
So we have verified that it is really enough to prove Proposition \ref{lm:main} 
which we will do in the next section after some preparation.
\end{proof}

\section{Main construction}

In this sections 
let the sets $K$, $L$, $M$,  the partitions $\{K_n:n\in {\omega}\}$
and $\{L_n:n\in {\omega}\}$, and the families $\mc K$, $\mc L$ and $\mc M$ be fixed as 
in the proof of  Theorem \ref{tm:mainmain}.

\begin{definition}\label{df:approx}
A pair ${\mathbb D}=\descris{}$ is a {\em description} iff 
$\mc B=\{B_i:i\in {\nu}\}$ is a clopen subbase of a 0-dimensional $T_2$ topology 
${\tau}$ on $X$.

We say that $\descris{}$ is {\em countable} iff $X\cup {\nu}$ is countable.

If $\mbb D^{\alpha}$ is a description, we assume that  
$\mbb D^{{\alpha}}=\descris{{{\alpha}}}$, $\mc B^{{\alpha}}=\{B^{\alpha}_i:i<{\nu}^{\alpha}\}$,
and $\mc B^{{\alpha}}$ is the base of the topology ${\tau}^{\alpha}$.

We write 
\begin{displaymath}
\descris{0}\preceq \descris{1} 
\end{displaymath}
iff 
\begin{enumerate}[(i)]
\item $X^0\subs X^1$ and ${\nu}^0\subs {\nu}^1$,
\item $B^1_i\cap X^0=B^0_i$ for each $i\in {\nu}^0$,
\item $B^0_i\subs B^0_j$ iff $B^1_i\subs B^1_j$ for each $i,j\in{\nu}_0$,
\item $B^0_i\cap B^0_j=\empt $ iff $B^1_i\cap B^1_j=\empt$ for each $i,j\in {\nu}_0$.
\end{enumerate}

If ${\beta}$ is a limit ordinal and $\vec {\mathbb D}=\<\mathbb D^{\alpha}:{\alpha}<{\beta}\>$ is a 
$\preceq$-increasing sequence of descriptions, then  we take
\begin{enumerate}[(a)]
\item $X^{\beta}=\bigcup_{{\alpha}<{\beta}}X^{\alpha}$ and 
${\nu}^{\beta}=\bigcup_{{\alpha}<{\beta}}{\nu}^{\alpha}$,
\item $B^{\beta}_i=\bigcup\{B^{\alpha}_i:i\in {\nu}^{\alpha}\}$ for $i\in {\nu}^{\beta}$,
\end{enumerate}
and write 
\begin{displaymath}
\lim \vec{\mathbb D}=\descri{{\beta}}.
\end{displaymath}
\end{definition}

The statement of the following lemma is straightforward.

\begin{lemma}\label{lm:descri-main}
   $\lim \vec{\mathbb D}$ is a description and ${\mathbb D}^{\alpha}\preceq \lim \vec{\mathbb D}
   $
for each ${\alpha}<{\beta}$.
\end{lemma}

\begin{definition}\label{df:appr}
   A triple $\mbb A=\<X,\mc B,\mc F\>$ is an {\em approximation}
   iff 
   \begin{enumerate}[(i)]
   \item $X=K\cup L\cup M\cup {\gamma}$ for some ${\gamma}<{\omega}_1$,
   \item $\<X,\mc B\>$ is a countable description,  $\mc B=\{B_i:i<{\nu}\}$ for some ${\nu}<{\omega}_1$,
   \item $\mc F=\{F_y:y\in {\gamma}\}$, and 
   $F_y\in \mc K\cup \mc L \cup \mc M$,
   \item $F_y\to_{{\tau}} y$ for each $y\in {\gamma}$.
   \item For each $B\in \mc B$, $|B\cap M|={\omega}$  and  
   $\forall^\infty  n\in {\omega}$ 
   $|B\cap K_n|=|B\cap L_n|={\omega}$. 
   \end{enumerate}

   Write 
   \begin{displaymath}
   desc(\mbb A)=\<X,\mc B\>.
   \end{displaymath}

   If ${\mathbb A}^{\nu}$ is an approximation, we assume that  
   ${\mathbb A}^{\nu}=\appros{{\nu}}$,  and $\mc F^{\nu}=\{F^{\nu}_y:y\in {\gamma}^{\nu}\}$.

   Assume that ${\mathbb A}^k=\appros{k}$ are approximations for $k<2$.
   Write ${\mathbb A}^0\preceq  {\mathbb A}^1$ iff 
   \begin{enumerate}[(1)]
   \item $desc({\mathbb A}^0)\preceq desc({\mathbb A}^1)$,
   \item if $F^0_y\in \mc K$, then $F^1_y=F^0_y$,   
   \item if $F^0_y\in \mc L\cup \mc M$, then $F^1_y\subs^*F^0_y$.
   \end{enumerate}
   \end{definition}
   
   \begin{lemma}\label{lm:approx}
      Assume that $\mathbb A=\<{\mathbb A}^{\alpha}:{\alpha}<{\beta}\>$ is a $\preceq$-increasing sequence of 
      approximations. 
      Write $\<X^{\beta},B^{\beta}\>=\lim \<desc({\mathbb A}^{\alpha}):{\alpha}<{\beta}\>$.
      
      \noindent (1)   
      For each ${\alpha}<{\beta}$ and $y\in {\gamma}_{\beta}$,
      \begin{enumerate}[(a)]
      \item   if $F^{\alpha}_y\in \mc K$, then $F_y\to_{{\tau}_{\beta}} y$.
      \item if $F^{\alpha}_y\in \mc L\cup \mc M$, then $y\in acc_{{\tau}_{\beta}}(F^{\alpha}_y) $.
      \end{enumerate}

      \noindent (2)   If  ${\beta}$ is countable, 
      there is a countable  approximation 
      ${\mathbb A}^{\beta}=\appros{{\beta}}$ such that 
      ${\mathbb A}^{\alpha}\prec {\mathbb A}^{\beta}$ for each ${\alpha}<{\beta}$
      (recall that $\<X^{\beta},B^{\beta}\>=\lim(\<\<X^{\alpha},B^{\alpha}\>:{\alpha}<{\beta}\>)$).
      \end{lemma}

   \begin{proof}[Proof of the Lemma]
   (1) is straightforward.
   To prove (2) we should define $F^{\beta}_y$ for $y\in {\gamma}^{\beta}$. 
   If $F^{\alpha}_y\subs K$ for some ${\alpha}<{\beta}$, then let 
   $F^{\beta}_y=F^{\alpha}_y$.

   If $F^{\alpha}_y\subs L\cup M$, then $\<F^{\alpha}_y:y\in {\gamma}_{\alpha}, {\alpha}<{\beta}\>$
   is a countable, mod-finite decreasing sequence of infinite sets. 
   So we can choose $F^{\beta}_y$ such that
   $F^y_{\beta}\subs^* F^{\alpha}_y$ for each ${\alpha}<{\beta}$ with 
   $y\in {\gamma}_{\alpha}$.
   
   Then $\mbb A^{\beta}=\<X^{\beta},\mc B^{\beta},\{F^{\beta}_y:y\in {\gamma}^{\beta}\}\>$
   is a suitable approximation.  
   \end{proof}

\begin{definition}\label{df:d}
If $\mbb A=\appros{}$ is an approximation and $n\in {\omega}$, 
let 
\begin{displaymath}
\dset(\mbb A,n)=\{B\cap K_j,B\cap L_j, B\cap M:B\in \mc B,n\le j<{\omega}\}.
\end{displaymath}
\end{definition}

\begin{lemma}\label{lm:a0a1partition}
If $\mbb A^0=\appros{0}$ is an approximation and 
$\<H_0,H_1\>$ is a partition of $X^0$ such that for some $n\in {\omega}$, 
\begin{enumerate}[(1)]
\item $|H_i\cap D|={\omega}$ for each $i<2$ and $D\in \dset(\mbb A^0,n)\cap {[X^0]}^{{\omega}}$,
\item $F^0_y\subs^* H_i$ for each $i<2$ and $y\in H_i\cap {\gamma}$,  
\end{enumerate}
then there is an approximation $\mbb A^1\ge \mbb A^0$ such that 
$H_0,H_1\in \mc B^1$.
\end{lemma}

\begin{proof}
   Let $\mc B_1=\{B^1_{\zeta}:{\zeta}<{\gamma}^1\}$  be an enumeration of 
   $\mc B^0\cup\{B\cap H_i:B\in \mc B^0,i<2\}\cup\{H_0,H_1\}$
   such that $B^1_{\zeta}=B^0_{\zeta}$ for ${\zeta}<{\gamma}^0$.
Then $\mc B^1$ is a base of a  
0-dimensional $T_2$  topology on $X$
and $H_0,H_1\in \mc B^1$.

Then $\mbb A^1=\<X^0,\mc B^1,\mc F^0\>$ is an approximation 
which  meets the requirements. 
\end{proof}

\begin{lemma}\label{lm:refine1}
   Assume that  $\mbb A^0=\appros{0}$ is an approximation, 
$n\in {\omega}$
$k\in K$, $\ell\in L$, $m\in M$, $y\in {\gamma}^0$,
$K'\in \mc K$, $L'\in \mc L$ such that 
$F_y\cap (K'\cup L')$ is finite.
Then there is an approximation $\mbb A^1\ge \mbb A^0$ such that $X^1=X^0$, $\mc F^1=\mc F^0$, and 
\begin{displaymath}
\{k,\ell,m,y\}\cap acc(K'\cup K_{\le n}\cup L'\cup L_{\le n}\cup {\gamma}^0,{\tau}^1)=\empt.
\end{displaymath}
\end{lemma}

\begin{proof}
Choose  $K''\in \mc K$ and $L''\in \mc L$ such that 
$K''\supset K'$, $L''\supset L'$ and 
$F^0_z\subs^* K''\cup L''$ provided $F^0_z\subs K\cup L$ for each $z\in {\gamma}^0$.

Choose $M''\in \mc M$ such that 
$F^0_z\subs M$ implies $F^0_z\subs ^* M''$ for each $z\in {\gamma}^0$, and 
$|B\setm M''|={\omega}$ for each $B\in \mc B$.

Write $N''=K''\cup L''\cup M''\cup K_{<n}\cup L_{<n}$.

Since $D\setm  N''$ is infinite for each 
$D\in \dset(\mbb A^0,n)$, there is a partition 
$\<G_0,G_1\>$ of $(K\cup L \cup M)\setm N''$
such that $G_i\cap D$ is infinite for each infinite 
$D\in \dset(\mbb A^0,n)$.

Let 
\begin{align*}
H_0=&G_0\cup F_y\cup \{y,k,\ell.m\},\\ 
H_1=&G_1\cup(N''\setm  (F_y\cup\{k,\ell,m\})\cup ({\gamma}^0\setm\{y\}).
\end{align*}
Clearly, $H_1=X^0\setm H_0$.
Then $\<H_0,H_1\>$ is a partition of $X^0$ and we can apply Lemma \ref{lm:a0a1partition}
to obtain $\mbb A^1$.
\end{proof}

\begin{lemma}\label{lm:refine2}
   Assume that  $\mbb A^0$ is a countable approximation, and $K'\in \mc K$.
   Then there is a countable approximation $\mbb A^1\ge \mbb A^0$ such that
   $F^1_y\cap K$ is infinite for some $y\in {\gamma}^1$. 
\end{lemma}

\begin{proof}[Proof of Lemma \ref{lm:refine2}]
We can assume that $K'\cap F^0_y$ is finite for each $y\in {\gamma}^0$.

Next, applying Lemma \ref{lm:refine1} ${\omega}$ times, 
we can get $\mbb A^2\ge \mbb A^0$ such that, $X^2=X^1$, $\mc F^2=\mc F^1$,  $K'$ is closed discrete in ${\tau}^2$
and the $K_n$-s and $L_n$-s are also closed discrete  
(in each step we guarantee that a point $x\in X^0$ is not an accumulation point of 
$K'\cup K_n\cup L_n$ ). 

Then, by induction on $n\in {\omega}$, we can pick $k_n\in K'$ and $x_n\in U_n\in \mc B^0$ such that 
$\{B_n:n<{\omega}\}$ is a locally finite family of disjoint open sets
such that
\begin{displaymath}
\forall B\in \mc B^2\ \big( \forall^\infty  n\in {\omega}\  U_n\subs B\ \lor
\forall^\infty  n\in {\omega}\  U_n\cap B=\empt\ \big).
\end{displaymath} 
Let ${\gamma}^1={\gamma}^2+1$, ${\nu}^1={\nu}^2+{\omega}$,
and let 
\begin{displaymath}
B^1_{{\nu}^2+n}=\{{\nu}^2\}\cup\bigcup\{U_j:j\ge n\}.
\end{displaymath}
For ${\zeta}<{\nu}^2$ let
\begin{displaymath}
B^1_{\zeta}
=\left\{\begin{array}{ll}
{B^2_{\zeta}\cup\{{\gamma}^0\}}&\text{if $ \forall^\infty  n\in {\omega}\ U_n\subs B^2_{\zeta}$,}\\\\
{B^2_{\zeta}}&\text{if $\forall^\infty  n\in {\omega}\ U_n\cap B^2_{\zeta}=\empt $}.
\end{array}\right.
\end{displaymath}
Let $F^1_{{\gamma}^0}=\{x_n:n<{\omega}\}$.
Then $\mbb A^1=\appros{1}$ meets the requirements. 
\end{proof}

\begin{lemma}\label{lm:refine3}
   Assume that  $\mbb A^0$ is a countable approximation, and $L'\in \mc L$.
   Then there is a countable approximation $\mbb A^1\ge \mbb A^0$ such that
   $F^1_y\subs^* L'$  for some $y\in {\gamma}^1$. 
\end{lemma}

\begin{proof}[Proof of Lemma \ref{lm:refine3}]
Imitating the  proof of Lemma \ref{lm:refine2}
we obtain  $\mbb A^2\ge \mbb A^0$ such that
$F^2_y\cap  L'$ is infinite  for some $y\in {\gamma}^2$.
Define $\mbb A^1=\< X^2,\mc B^2,\mc F^1\>$
such that 
\begin{displaymath}
F^1_y=\left\{\begin{array}{ll}
{F^2_y\cap L'}&\text{\text {if ${F^2_y\cap L'}$ is infinite},}\\\\
{F^2_y}&\text{\text {if {$F^2_y\cap L'$} is finite.}}
\end{array}\right.
\end{displaymath}
Then $\mbb A^1$ satisfies the requirements. 
\end{proof}

\begin{lemma}\label{lm:refine4pre}
   Assume that  $\mbb A^0$ is a countable approximation, and $M'\in \mc M$.
   Then there is a countable approximation $\mbb A^1\ge \mbb A^0$ such that
   $M'$ contains an infinite closed discrete subset.  
\end{lemma}

\begin{proof}[Proof of Lemma \ref{lm:refine4}]
   We can assume that $M'\to_{{\tau}^0}x$ for some $x\in X$.
   Pick $K'\in \mc K $ such that $F^0_y\subs^* K'$ for each 
   $y\in {\gamma}^0$ with $F^0_y\subs K$. 
   Let 
   \begin{displaymath}
      \mc D= \dset(\mbb A^0,0)\cup(\{B\cap F^0_z:z\in {\gamma}^0, B\in \mc B^0\}\cap {[L\cup M]}^{{\omega}})\cup\{M'\cap B:B\in \mc B^0\}.
      \end{displaymath}

      Let $\<G_0,G_1\>$ be a partition of $(K\cup L\cup M)\setm K'$ such that 
   $G_i\cap D$ is infinite for each infinite $D\in \mc D$. 
   
   Let $H_0=G_0\cup K'\cup {\gamma}^0$ and $H_1=X^0\setm H_0$.
   
   Let $\mc B_1=\{B^1_{\zeta}:{\zeta}<{\gamma}^1\}$  be an enumeration of 
   $\mc B^0\cup\{B\cap H_i:B\in \mc B^0,i<2\}\cup\{H_0,H_1\}$
   such that $B^1_{\zeta}=B^0_{\zeta}$ for ${\zeta}<{\gamma}^0$.
   Then $\mc B^1$ is a base of a  
   0-dimensional $T_2$  topology on $X$
   and $H_0,H_1\in \mc B^1$.
   
   Then $\mbb A^1=\<X^0,\mc B^1,\mc F^0\>$ is an approximation 
   which  meets the requirements. 
because $M'\cap H_1$ is closed discrete. 
   \end{proof}

\begin{lemma}\label{lm:refine4}
   Assume that  $\mbb A^0$ is a countable approximation, and $M'\in \mc M$.
   Then there is a countable approximation $\mbb A^1\ge \mbb A^0$ such that
   $F^1_y\subs^* M'$  for some $y\in {\gamma}^1$. 
\end{lemma}

\begin{proof}[Proof of Lemma \ref{lm:refine4}]

   By Lemma \ref{lm:refine4pre} we can assume that $M'$ is closed discrete. 

   Then, by induction on $n\in {\omega}$, we can pick $k_n\in M'$ and $x_n\in B_n\in \mc B^0$ such that 
$\{B_n:n<{\omega}\}$ is a locally finite family of disjoint open sets
such that
\begin{displaymath}
\forall B\in \mc B^2\ \big( \forall^\infty  n\in {\omega}\  U_n\subs B\ \lor
\forall^\infty  n\in {\omega}\  U_n\cap B=\empt\ \big).
\end{displaymath} 
Let ${\gamma}^1={\gamma}^2+1$, ${\nu}^1={\nu}^2+{\omega}$,
and let 
\begin{displaymath}
B^1_{{\nu}^2+n}=\{{\nu}^2\}\cup\bigcup\{U_j:j\ge n\}.
\end{displaymath}
For ${\zeta}<{\nu}^2$ let
\begin{displaymath}
B^1_{\zeta}
=\left\{\begin{array}{ll}
{B^2_{\zeta}\cup\{{\gamma}^0\}}&\text{if $ \forall^\infty  n\in {\omega}\ U_n\subs B^2_{\zeta}$,}\\\\
{B^2_{\zeta}}&\text{if $\forall^\infty  n\in {\omega}\ U_n\cap B^2_{\zeta}=\empt $}.
\end{array}\right.
\end{displaymath}
Let $F^2_{{\gamma}^0}=\{x_n:n<{\omega}\}$.
Then $\mbb A^2=\appros{2}\le \mbb A^1$ and there is $y\in {\gamma}^2$
such that $M'\cap F^2_y$ is infinite.

Define $\mbb A^1=\< X^2,\mc B^2,\mc F^1\>$
such that 
\begin{displaymath}
F^1_y=\left\{\begin{array}{ll}
{F^2_y\cap L'}&\text{\text {if ${F^2_y\cap L'}$ is infinite}}\\\\
{F^2_y}&\text{\text {if {$F^2_y\cap L'$} is finite.}}
\end{array}\right.
\end{displaymath}
Then $\mbb A^1$ satisfies the requirements. 

\end{proof}

\begin{lemma}\label{lm:refine5}
   Assume that  $\mbb A^0$ is a countable approximation, 
   $y\in {\gamma}^0$
   and $L'\in \mc L$.
   Then there is a countable approximation $\mbb A^1\le \mbb A^0$ such that
   $'\not\to_{{\tau}^1}y$. 
\end{lemma}

\begin{proof}[Proof of Lemma \ref{lm:refine5}]
Let $K'\in \mc K$ such that $F^0_z\subs^* K'$ for each $z\in {\gamma}^0$ with $F^0_z\subs K$.

Let 
\begin{displaymath}
   \mc D= \dset(\mbb A^0,n)\cup(\{B\cap F^0_z:z\in {\gamma}^0, B\in \mc B^0\}\cap {[L\cup M]}^{{\omega}})\cup\{L'\cap B:B\in \mc B^0\}.
   \end{displaymath}

   Let $\<G_0,G_1\>$ be a partition of $(K\cup L\cup M)\setm K_0$ such that 
$H_i\cap D$ is infinite for each infinite $D\in \mc D$. 

Let $H_0=G_0\cup K'\cup Y$ and $H_1=X^0\setm H_0$.

Let $\mc B_1=\{B^1_{\zeta}:{\zeta}<{\gamma}^1\}$  be an enumeration of 
$\mc B^0\cup\{B\cap H_i:B\in \mc B^0,i<2\}\cup\{H_0,H_1\}$
such that $B^1_{\zeta}=B^0_{\zeta}$ for ${\zeta}<{\gamma}^0$.
Then $\mc B^1$ is a base of a  
0-dimensional $T_2$  topology on $X$
and $H_0,H_1\in \mc B^1$.

Then $\mbb A^1=\<X^0,\mc B^1,\mc F^0\>$ is an approximation 
which  meets the requirements.

\end{proof}

\begin{proof}[Proof of Proposition \ref{lm:main}]
   Let 
   \begin{multline*}
   Task=\{\<1,Y\>:Y\in \mc K\cup \mc L\cup \mc M\}\cup \\\{\<2,Z\>:Z\in \mc K\cup \mc L\cup \mc M\cup \mc {\omega}_1\cup\{L_n,K_n:n<{\omega}\}\}\\\cup
   \{\<3,T\>:T\in \mc L\cup \mc M\}.
   \end{multline*}
and fix an ${\omega}_1$-abundant enumeration $\{t_{\zeta}:{\zeta}<{\omega}_1\}$   of 
   $Task$.

   We will define a $\preceq$-increasing sequence $\<\mbb A^{\alpha}=\appros{{\alpha}}:{\alpha}<{\omega}_1\>$
 of 
countable approximations by transfinite recursion. 

Let $\mbb A^0$ be a countable approximation such that ${\gamma}^0=0$.

If ${\zeta}$ is a limit ordinal, apply Lemma \ref{lm:approx}(2) to obtain 
$\mbb A^{\zeta}$ from $\<\mbb A^{\eta}:{\eta}<{\zeta}\>$.

Assume that ${\zeta}={\eta}+1$ and we have $\mbb A^{\eta}$.

\noindent{\bf Case 1.} $t_{\eta}=\<1,Y\>$. 

If $Y\in \mc K$, then apply Lemma \ref{lm:refine2} to find
$\mbb A^{\zeta}$ such that $Y\cap F^{\zeta}_y$ is infinite for some 
$y<{\gamma}^{\zeta}$. 

If $Y\in \mc L\cup \mc M$, then apply Lemma \ref{lm:refine3} or Lemma \ref{lm:refine5}to find
$\mbb A^{\zeta}$ such that $F^{\zeta}_y\subs F$  for some 
$y<{\gamma}^{\zeta}$.

\noindent{\bf Case 2.} $t_{\eta}=\<2,Z\>$.

If $Z\in \mc K\cup \mc L$, then apply Lemma \ref{lm:refine1} ${\omega}$ times 
to obtain $\mbb A^{\zeta}$ such that $acc(Z,{\tau}^{\zeta})\subs {\gamma}^{\eta}$.

If $Z\in {\omega}_1\cup \{K_n,L_n:n<{\omega}\}$, then apply Lemma \ref{lm:refine1} ${\omega}$ times 
to obtain $\mbb A^{\zeta}$ such that $acc(Z,{\tau}_{\zeta})=\empt$.

\noindent{\bf Case 3.} $t_{\eta}=\<3,T\>$.
Apply Lemma \ref{lm:refine4pre} or Lemma \ref{lm:refine5}
to obtain $\mbb A^{\zeta}$ such that $T$ does not converge in ${\tau}^{\zeta}$.

FInally let $\<X,\mc B \>=\lim\<\<X^{\zeta},\mc B^{\zeta}\>:{\zeta}<{\omega}_1\>$.
Then the space $\<X^{{\omega}_1},{\tau}^{{\omega}_1}\>$
satisfies the requirements.

\end{proof}

\end{document}